\documentclass[10pt]{article}

\usepackage{latexsym}
\usepackage{amssymb}
\usepackage{amsfonts}

\newtheorem{thm}{Theorem}
\newtheorem{lem}{Lemma}
\newtheorem{cor}{Corollary}
\newtheorem{rem}{Remark}
\newtheorem{prop}{Proposition}
\newtheorem{prob}{Problem}
\newtheorem{defn}{Definition}

\newcommand{\bpr}{\noindent{\bf PROOF\/. }}
\newcommand{\epr}{\hspace*{\fill}$\square$\medskip}

\newcommand{\bprA}{\noindent{\bf PROOF OF THEOREM 3\/. }}

\author{Houmem Belkhechine\\
Facult\'e des Sciences de Gab\`es \\Cit\'e Riadh, Zirig \\6072 Gab\`es \\Tunisie\\
{\tt houmem@gmail.com}\and
Imed Boudabbous\\
Institut Pr\'eparatoire aux \'Etudes d'Ing\'enieurs de Sfax \\
Route Menzel Chaker Km 0.5 \\
3018 Sfax \\Tunisie\\
{\tt imed.boudabbous@gmail.com}}

 \title{Indecomposable tournaments \\and their indecomposable subtournaments\\ on 5 and 7 vertices}

\date{}

\begin{document}

\maketitle

\begin{abstract}
 Given a tournament $T=(V,A)$, a subset
$X$ of $V$ is an interval of $T$ provided that for every $a, b\in X$
and $ x\in V-X$, $(a,x)\in A$ if and only if $(b,x)\in A$. For
example, $\emptyset$, $\{x\}(x\in V)$ and $V$ are intervals of $T$,
called trivial intervals. A tournament, all the intervals of which
are trivial, is indecomposable; otherwise, it is decomposable. A
critical tournament is an indecomposable tournament $T$ of
cardinality $\geq 5$ such that for any vertex $x$ of $T$, the
tournament $T-x$ is decomposable. The critical tournaments are of
odd cardinality and for all $n \geq 2$ there are exactly three
critical tournaments on $2n+1$ vertices denoted by $T_{2n+1}$,
$U_{2n+1}$ and $W_{2n+1}$. The tournaments $T_{5}$, $U_{5}$ and
$W_{5}$ are the unique indecomposable tournaments on 5 vertices. We
say that a tournament $T$ embeds into a tournament $T'$ when $T$ is
isomorphic to a subtournament of $T'$. A diamond is a tournament on
4 vertices admitting only one interval of cardinality 3. We prove
the following theorem: if a diamond and $T_{5}$ embed into an
indecomposable tournament $T$, then $W_{5}$ and $U_{5}$ embed into
$T$. To conclude, we prove the following: given an indecomposable
tournament $T$, with $\mid\!V(T)\!\mid \geq 7$, $T$ is critical if
and only if the indecomposable subtournaments on 7 vertices of $T$ are isomorphic to 
one and only one of the tournaments $T_{7}$, $U_{7}$ and  $W_{7}$.
\end{abstract}

\noindent {\bf Key words:} Tournament; Indecomposable; Critical; Embedding.

\section{Basic definitions}

A {\it tournament} $T=(V(T),A(T))$ or $(V,A)$ consists of a finite
{\it vertex} set $V$ with an {\it arc} set $A$ of ordered pairs of
distinct vertices satisfying: for $x, y \in V$, with $x \neq y$,
$(x, y) \in A$ if and only if $(y, x) \notin A$. The {\it
cardinality} of T is that of $V(T)$ denoted by $\mid\!V(T)\!\mid$.
For two distinct vertices $x$ and $y$ of a tournament $T$, $x
\longrightarrow y$ means that $(x,y) \in A(T)$. For $x \in V(T)$ and
$Y \subset V(T)$, $x \longrightarrow Y$ (rep. $Y \longrightarrow x$)
signifies that for every $y \in Y$, $x \longrightarrow y$ (resp. $y
\longrightarrow x$).  Given a vertex $x$  of a tournament $T = (V,
A)$, $N_{T}^{+}(x)$   denotes  the set $\{y \in V: x \longrightarrow
y \}$. The {\it score} of $x$ (in $T$), denoted by $s_{T}(x)$, is
the cardinality of $N_{T}^{+}(x)$. A tournament is {\it regular} if
all its vertices share the same score. A {\it transitive} tournament
or {\it total order} is a tournament $T$ such that for $x$, $y$,
$z\in V(T)$, if $x \longrightarrow y$ and $y \longrightarrow z$,
then $x \longrightarrow z$. For two distinct vertices $x$ and $y$ of
a total order $T$, $x < y$ means that $x \longrightarrow y$. We
write $T = a_{0} < \cdots < a_{n}$ to mean that $T$ is the total
order defined on $V(T) = \{a_{0}, \ldots, a_{n}\}$ by $A(T) =
\{(a_{i}, a_{j}) : i < j\}$.

The notions of isomorphism, of subtournament and of embedding are
defined in the following manner. First, let $T = (V, A)$ and $T' =
(V', A')$ be two tournaments. A one-to-one correspondence $f$ from
$V$ onto $V'$ is an {\it isomorphism} from $T$ onto $T'$ provided
that for $x, y \in V$, $(x, y) \in A$ if and only if $(f(x), f(y))
\in A'$. The tournaments $T$ and $T'$ are then said to be {\it
isomorphic}, which is denoted by $T \simeq T'$. Moreover, an
isomorphism from a tournament $T$ onto itself is called an {\it
automorphism} of $T$. The automorphisms of $T$ form a subgroup of
the permutation group of $V(T)$, called the {\it automorphism group}
of $T$. Second, given a tournament $T = (V, A)$, with each subset
$X$ of $V$ is associated the {\it subtournament} $T(X) = (X, A \cap
(X \times X))$ of $T$ {\it induced} by $X$. For $x \in V $, the
subtournament $T(V - \{x\})$ is denoted by $T-x$. For tournaments
$T$ and $T'$, if $T'$ is isomorphic to a subtournament of $T$, then
we say that $T'$ {\it embeds} into $T$. Otherwise, we say that $T$
{\it omits} $T'$. The {\it dual} of a tournament $T=(V, A)$ is the
tournament obtained from $T$ by reversing all its arcs. This
tournament is denoted by $T^{\star}=(V, A^{\star})$, where
$A^{\star} = \{(x, y): \ (y,x) \in A\}$. A tournament $T$ is then
said to be {\it self-dual} if $T$ and $T^{\star}$ are isomorphic.

The indecomposability plays an important role in this paper. Given a
tournament $T = (V, A)$, a subset $I$ of $V$ is an \emph{interval}
(\cite{F1}, \cite{I}, \cite{ST}) (or a {\it clan} \cite{E} or an
{\it homogeneous subset} \cite{G}) of $T$ provided that for every $x
\in V-I$, $x \rightarrow I$ or $I \rightarrow x$. This definition
generalizes the notion of interval of a total order. Given a
tournament $T = (V, A)$, $\emptyset $, $V$ and $\{x\}$, where $x \in
V$, are clearly intervals of $T$, called {\it trivial} intervals. A
tournament is then said to be {\it indecomposable} (\cite{I},
\cite{ST}) (or {\it primitive} \cite{E}) if all of its intervals are
trivial, and is said to be {\it decomposable} otherwise. For
instance, the 3-cycle $C_{3}  = (\{0, 1, 2\}, \{(0, 1), (1, 2), (2,
0)\})$ is indecomposable whereas a total order of cardinality $\geq
3$  is decomposable. Let us mention the following relationship
between indecomposability and duality. The tournaments $T$ and
$T^{\star}$ have the same intervals and, thus, $T$ is indecomposable
if and only if $T^{\star}$ is indecomposable.

\section{The critical tournaments}

An indecomposable tournament $T = (V, A)$ is said to be {\it
critical} if $\mid\!V\!\mid > 1$ and for all $x \in V$, $T-x$ is
decomposable. In order to present our main results and to present
the characterization of the critical tournaments due to J.H.~
Schmerl and W.T. Trotter \cite{ST}, we introduce the tournaments
$T_{2n + 1}$, $U_{2n + 1}$ and $W_{2n + 1}$ defined on $2n + 1$
vertices, where $n\geq2$, as follows:

\begin{itemize}
\item The tournament $T_{2n+1}$ is the tournament defined on
$\mathbb{Z}/(2n+1)\mathbb{Z}$ by $A(T_{2n+1}) = \{(i,j): j-i \in
\{1, \ldots, n\}\}$, so that, $T_{2n+1}(\{0,\ldots, n\}) = 0< \cdots
<n$, $T_{2n+1}(\{n+1, \ldots, 2n\}) = n+1< \cdots < 2n$ and for $i
\in \{0, \ldots, n-1\},\ \{i+1, \ldots, n\} \longrightarrow i+n+1
\longrightarrow\{0, \ldots, i\}$ (see Figure 1).

\item The tournament $U_{2n+1}$ is obtained from $T_{2n+1}$ by reversing the arcs of $T_{2n+1}(\{n+1,
\ldots, 2n\})$. Therefore, $U_{2n+1}$ is defined on $\{0,\ldots,
2n\}$ as follows: $U_{2n+1}(\{0,\ldots, n\}) = 0< \cdots <n$,
$U_{2n+1}^{\star}(\{n+1, \ldots, 2n\}) = n+1< \cdots < 2n$ and for
$i \in \{0, \ldots, n-1\},\ \{i+1, \ldots, n\} \longrightarrow i+n+1
\longrightarrow\{0, \ldots, i\}$ (see Figure 2).

\item The tournament $W_{2n+1}$ is defined on $\{0,\ldots, 2n\}$ in the following manner: $W_{2n+1}-2n = 0< \cdots <2n-1$
and  $\{1, 3, \ldots, 2n-1\} \longrightarrow 2n \longrightarrow \{0,
2, \ldots, 2n-2 \}$ (see Figure 3).
\end{itemize}

\begin{figure}[ht]
\begin{center}
\setlength{\unitlength}{1cm}
\begin{picture}(11,5)
\put(0,1){0} \put(1.5,1){1} \put(3.5,1){i} \put(2,1.1){.}
\put(2.5,1.1){.} \put(3,1.1){.} \put(5,1){i+1} \put(5.9,1.1){.}
\put(6.2,1.1){.} \put(6.5,1.1){.} \put(7,1){$n-1$}
\put(9.2,1){$n$} \put(1.6,0.8){\line(2,-1){0.5}}
\put(2.1,0.55){\vector(1,0){3.2}} \put(0.1,0.8){\line(3,-2){1.14}}
\put(1.24,0.04){\vector(1,0){6.23}} \put(0.3,1.1){\vector(1,0){1}}
\put(3.8,1.1){\vector(1,0){1}} \put(8,1.1){\vector(1,0){1}}
\put(0.3,2.5){n+1} \put(3.8,2.5){i+n+1} \put(8.3,2.5){2n}
\put(1.9,2.6){\vector(1,0){1}} \put(1.15,2.6){.} \put(1.4,2.6){.}
\put(1.65,2.6){.} \put(3.1,2.6){.} \put(3.32,2.6){.}
\put(3.54,2.6){.}
\put(5.8,2.6){\vector(1,0){1.5}} \put(5,2.6){.} \put(5.25,2.6){.}
\put(5.5,2.6){.} \put(7.55,2.6){.} \put(7.8,2.6){.}
\put(8.05,2.6){.}
\put(4.2,2.9){\line(2,1){0.52}} \put(4.7,3.15){\vector(1,0){1.77}}
\put(0.65,2.9){\line(2,1){1.32}}
\put(1.95,3.55){\vector(1,0){1.5}}
\put(0.3,2.3){\vector(-1,-4){0.23}}
\put(1.55,1.4){\vector(-1,1){0.9}}
\put(2.3,1.3){\vector(-3,2){1.4}}
\put(3.8,2.3){\vector(-1,-1){1}} \put(4,2.3){\vector(-1,-3){0.3}}
\put(5.25,1.4){\vector(-1,1){0.9}}
\put(6.1,1.3){\vector(-3,2){1.5}}
\put(8.45,2.3){\vector(-1,-1){1}}
\put(8.2,2.3){\vector(-3,-2){1.5}}
\put(9.57,1.4){\vector(-1,1){0.9}}
\end{picture}
\end{center}
\caption{$T_{2n+1}$.}
\end{figure}

\begin{figure}[ht]
\begin{center}
\setlength{\unitlength}{1cm}
\begin{picture}(11,5)
\put(0,1){0} \put(1.5,1){1} \put(3.5,1){i} \put(2,1.1){.}
\put(2.5,1.1){.} \put(3,1.1){.} \put(5,1){i+1} \put(5.9,1.1){.}
\put(6.2,1.1){.} \put(6.5,1.1){.} \put(7,1){$n-1$}
\put(9.2,1){$n$} \put(1.6,0.8){\line(2,-1){0.5}}
\put(2.1,0.55){\vector(1,0){3.2}} \put(0.1,0.8){\line(3,-2){1.14}}
\put(1.24,0.04){\vector(1,0){6.23}} \put(0.3,1.1){\vector(1,0){1}}
\put(3.8,1.1){\vector(1,0){1}} \put(8,1.1){\vector(1,0){1}}
\put(0.3,2.5){n+1} \put(3.8,2.5){i+n+1} \put(8.3,2.5){2n}
\put(2.9,2.6){\vector(-1,0){1}} \put(1.15,2.6){.} \put(1.4,2.6){.}
\put(1.65,2.6){.} \put(3.1,2.6){.} \put(3.32,2.6){.}
\put(3.54,2.6){.}
\put(7.3,2.6){\vector(-1,0){1.5}} \put(5,2.6){.} \put(5.25,2.6){.}
\put(5.5,2.6){.} \put(7.55,2.6){.} \put(7.8,2.6){.}
\put(8.05,2.6){.}
\put(4.2,2.9){\line(-2,1){0.5}} \put(3.7,3.15){\vector(-1,0){1.3}}
\put(8.5,2.9){\line(-2,1){1.32}}
\put(7.2,3.55){\vector(-1,0){1.9}}
\put(0.3,2.3){\vector(-1,-4){0.23}}
\put(1.55,1.4){\vector(-1,1){0.9}}
\put(2.3,1.3){\vector(-3,2){1.4}}
\put(3.8,2.3){\vector(-1,-1){1}} \put(4,2.3){\vector(-1,-3){0.3}}
\put(5.25,1.4){\vector(-1,1){0.9}}
\put(6.1,1.3){\vector(-3,2){1.5}}
\put(8.45,2.3){\vector(-1,-1){1}}
\put(8.2,2.3){\vector(-3,-2){1.5}}
\put(9.57,1.4){\vector(-1,1){0.9}}
\end{picture}
\end{center}
\caption{$U_{2n+1}$.}
\end{figure}


\begin{figure}[ht]
\begin{center}
\setlength{\unitlength}{1cm}
\begin{picture}(11,5)
\put(0,1){0} \put(1.5,1){1} \put(3.5,1){2i} \put(2,1.1){.}
\put(2.5,1.1){.} \put(3,1.1){.} \put(5,1){2i+1} \put(6,1.1){.}
\put(6.3,1.1){.} \put(6.6,1.1){.} \put(7,1){$2n-2$}
\put(9.2,1){$2n-1$} \put(1.6,0.8){\line(2,-1){0.5}}
\put(2.1,0.55){\vector(1,0){3.2}} \put(0.1,0.8){\line(3,-2){1.14}}
\put(1.24,0.04){\vector(1,0){6.23}} \put(0.3,1.1){\vector(1,0){1}}
\put(3.9,1.1){\vector(1,0){1}} \put(8.1,1.1){\vector(1,0){1}}
\put(4.2,3){2n}   
\put(3.8,3){\vector(-2,-1){3.4}}  
\put(1.9,1.4){\vector(3,2){2.1}}  
\put(4.25,2.8){\vector(-1,-3){0.475}}  
\put(5.35,1.4){\vector(-2,3){0.95}}  
\put(4.7,2.8){\vector(2,-1){2.9}}   
\put(9.7,1.4){\vector(-3,1){4.8}}   
\end{picture}
\end{center}
\caption{$W_{2n+1}$.}
\end{figure}

\begin{thm}[\cite{ST}]\label{t1} Up to isomorphism, the critical tournaments of cardinality $\geq 5$ are
the tournaments $T_{2n+1}$, $U_{2n+1}$ and $W_{2n+1}$, where $n
\geq 2$.
\end{thm}

Notice that the critical tournaments are self-dual.

\section{The tournaments $T_{5}$, $U_{5}$ and $W_{5}$ in an indecomposable tournament}
We study  the indecomposable tournaments according to their
indecomposable subtournaments on 5 vertices. A recent result on our
topic is a characterization of the indecomposable tournaments
omitting $W_{5}$ obtained  by B.J. Latka \cite{BJL}. In order to
recall this characterization, we introduce the {\it Paley}
tournament $P_{7}$ defined on $\mathbb{Z}/7\mathbb{Z}$ by $A(P_{7})
=\{(i, j): j-i \in \{1, 2, 4\}\}$. Notice that the tournaments
obtained from $P_{7}$ by deleting one vertex are isomorphic and
denote $P_{7} - 6$ by  $B_{6}$.

\begin{thm}[\cite{BJL}]\label{t2} Given a tournament $T$
of cardinality $\geq 5$, $T$ is indecomposable and omits $W_{5}$ if
and only if $ T $ is isomorphic to an element of $\{ B_{6},
P_{7}\}\cup \{T_{2n+1}: n \geq 2\}\cup \{U_{2n+1}: n \geq 2\}$.
\end{thm}

A {\it diamond} is a tournament on 4 vertices admitting only one
interval of cardinality 3. Up to isomorphism, there are exactly two
diamonds $D_{4}$ and $D_{4}^{\star}$, where $D_{4}$ is the
tournament defined on $\{0, 1, 2, 3\}$ by $D_{4}(\{0, 1,2\}) =
C_{3}$ and  $3\longrightarrow \{0, 1,2\}$.

The following theorem is the main result. This theorem is presented
in \cite{HI} without a detailed  proof.

\begin{thm}\label{t3} Given an indecomposable tournament $T$, if a diamond and $T_{5}$ embed into $T$, then $U_{5}$ and $W_{5}$ embed into $T$.
\end{thm}

 C. Gnanvo and P. Ille
\cite{GI} and G. Lopez and C. Rauzy \cite{LR} characterized the
tournaments omitting diamonds. In the indecomposable case they
obtained the following characterization.

\begin{prop}[\cite{GI,LR}]\label{p4}
Given an indecomposable tournament $T$ of cardinality $\geq 5$,
$T$ omits the diamonds $D_{4}$ and $D_{4}^{\star}$ if and only if
$T$ is isomorphic to $T_{2n+1}$ for some $n \geq 2$.
\end{prop}

\section{Proof of Theorem~\ref{t3}}

Before proving Theorem~\ref{t3}, we  introduce some notations and
definitions.

 \begin{defn}\label{d5} Given a tournament $T=(V,A)$, with each subset $X$ of $V$,
 such that $\mid\!X\!\mid \geq 3$ and $T(X)$ is indecomposable, are
 associated the following subsets of $V-X$.

 \begin{itemize}

\item $Ext(X) = \{x \in V-X: \ T(X \cup \{x\})$ is indecomposable$\}$.
\item $[X] = \{x \in V-X: \ x \rightarrow X$ or $X \rightarrow x \}$.

\item For every $u \in X$, $X(u) = \{x \in V-X: \ \{u, x\}$ is an interval
of $T(X \cup \{x\})\}$.

\end{itemize} \end{defn}

 \begin{lem}[\cite{E}]\label{l6} Let $T=(V,A)$ be a tournament and let $X$ be a subset
of $V$ such that $\mid\!X\!\mid \geq 3$ and $T(X)$ is
indecomposable.

\begin{enumerate}

\item The family $\{X(u): u \in X\} \cup \{Ext(X), [X]\}$
constitutes a partition of $V-X$.
\item Given $u\in X$, for all $x\in X(u)$ and for all $y\in V-(X\cup
X(u))$, if $T(X\cup \{x, y\})$ is decomposable, then $\{u, x\}$ is
an interval of $T(X\cup\{x, y\})$.
\item For every $x\in [X]$ and for every $y\in V-(X\cup [X])$, if $T(X\cup\{x,
y\})$ is decomposable, then $X\cup \{y\}$ is an interval of
$T(X\cup\{x, y\})$.
\item Given $x, y \in Ext(X)$, with $x \neq y$, if $T(X\cup\{x,
y\})$ is decomposable, then $\{x, y\}$ is an interval of
$T(X\cup\{x, y\})$. \end{enumerate}

\end{lem}

 The below result follows from Lemma~\ref{l6}.

\begin{prop}[\cite{E}]\label{p7} Let $T=(V,A)$ be an indecomposable tournament. If
$X$ is a subset of $V$, such that $\mid\!X\!\mid \geq 3$,
$\mid\!V-X\!\mid \geq 2$ and $T(X)$ is indecomposable, then there
are distinct  elements $x$ and $y$ of  $V-X$ such that $T(X\cup\{x,
y\})$ is indecomposable.

\end{prop}

\begin{cor}\label{c8}
Let $T = (V, A)$ be an indecomposable tournament such that $\mid V
\mid$ is even and $\mid V \mid \geq 6$. For each $x \in V$, there
is $y \in V - \{x\}$ such that $T - y$ is indecomposable.

\end{cor}
\bpr
As $T$ is indecomposable, there is $X \subset V$ such that $x \in X$
and $T(X) \simeq C_{3}$. Otherwise, $N_{T}^{+}(x)$ or $V- (\{x\}\cup
N_{T}^{+}(x))$ would be  non trivial intervals of $T$. Since $\mid V
\mid$ is even, by applying several times Proposition~\ref{p7} from the
indecomposable subtournament $T(X)$, we get a vertex $y \in V - X$
such that $T - y$ is indecomposable.
\epr

The 3-cycle $C_{3}$ is indecomposable and embeds into any
indecomposable tournament of cardinality $\geq3$ as observed in the
preceding proof. It follows, by Proposition~\ref{p7}, that any
indecomposable tournament $T$ of cardinality $\geq5$, admits an
indecomposable subtournament on 5 vertices. The indecomposable
tournaments on 5 vertices are critical because the four tournaments
on 4 vertices are decomposable. So let us mention the following
facts.

\begin{rem}\label{r9}\

\begin{itemize}

\item The indecomposable tournaments on 5 vertices are, up to
isomorphism, the three critical tournaments $T_{5}$, $U_{5}$ and
$W_{5} $.

\item There is no indecomposable tournament of cardinality $\geq5$
omitting each of the tournaments $T_{5}$, $U_{5}$ and $W_{5} $.
\end{itemize}
\end{rem}

The tournaments $T_{2n+1}$ play an important role in the proof of
Theorem~\ref{t3}. We recall some of their properties.

\begin{rem}\label{r10}\
\begin{itemize}

\item The tournaments $T_{2n+1}$ are regular: for all $i \in \{0, \ldots,
2n\}$, $s_{T_{2n+1}}(i) = n$;
\item For $0 \leq i \leq 2n$, the unique non trivial interval of
$T_{2n+1} -i$ is \{i+n, i+n+1\};

\item The automorphism group of $T_{2n+1}$ is generated by the permutation $\sigma: i \mapsto i+1$;

\item The permutation $\pi: i \mapsto -i$, is an isomorphism from $T_{2n+1}$ onto
its dual.

\end{itemize}

\end{rem}

Now we are ready to prove Theorem~\ref{t3}.

~

\bprA
Let $T = (V,A)$ be an indecomposable
tournament into which
  a diamond and $T_{5}$ embed. Consider a
minimal subset $X$ of $V$ such that $T(X)$ is indecomposable and a
diamond and $T_{5}$ embed into $T(X)$.  Now, let $Y$ be a maximal
subset of $X$ such that $T(Y) \simeq T_{2n+1}$ for some $n \geq 2$.
 We  establish that $\mid X \mid = 6$ by using the following
 observation. Consider a subset $Z$ of $X$ such that $T(Z) \simeq
 T_{2n+1}$ and assume that $Ext(Z) \cap X \neq \emptyset$. Let $x\in Ext(Z) \cap X
 $. We have $T(Z\cup\{x\})$  is indecomposable. Furthermore,
 as  $\mid Z\cup\{x\} \mid$ is even, a diamond embeds into
 $T(Z\cup\{x\})$ by Proposition~\ref{p4}. Since $T_{5}$ embeds into
 $T(Z\cup\{x\})$ as well, it follows from the minimality of $X$ that
 $X= Z\cup\{x\}$. As an immediate consequence, we have: if $Z$ is a
 subset of $X$ such that $T(Z) \simeq T_{2n+1}$ and $\mid X- Z \mid \geq 2$,
 then $Ext(Z) \cap X = \emptyset$. By Lemma~\ref{l6}, for every $x\in X-Z$,
 either $x\in [Z]$ or there is $u\in Z$ such that $x\in Z(u)$.

For a contradiction, suppose that $Ext(Y) \cap X = \emptyset$. By
Proposition~\ref{p7}, there are $x \neq y \in X-Y$ such that $T(Y \cup \{x,
y \})$ is indecomposable. Clearly, if $\{x, y\}\subseteq [Y]$, then
$Y$ would be a non trivial interval of $T(Y \cup \{x, y \})$. For
instance, assume that there is $v\in Y$ such that $y\in Y(v)$. By
Lemma~\ref{l6}, either there is $u\in Y$ such that $x\in Y(u)$ or $x\in
[Y]$. In each of  both instances, we obtain a contradiction.

First, suppose that there is $u \in Y$ such that $x \in Y(u)$. We
have $u\neq v$, otherwise $\{u, x, y\}$ would be a non trivial
interval of $T(Y \cup \{x, y \})$. By Remark~\ref{r10}, the automorphism
group of $T_{2n+1}$ is generated by $\sigma : i\mapsto i+1$.
Therefore, by interchanging $x$ and $y$, we can denote the element
of $Y$ by $0,\ldots, 2n$ in such a way that $T(Y)= T_{2n+1}$, $u =
0$ and $1\leq v\leq n$. Since $T(Y \cup \{x, y \})$ is
indecomposable and $0\longrightarrow v$, we get $y\longrightarrow x$
by Lemma~\ref{l6}. Consider $Z = (Y-\{0\}) \cup \{x\}$. We have $T(Z)
\simeq T_{2n+1}$ and, by the preceding observation, either $y\in
[Z]$ or there is $w \in Z$ such that $y \in Z(w)$. The first
instance is not possible because $\{v-2, v-1\} \cap Z \neq
\emptyset$ and $\{v-2, v-1\} \longrightarrow y \longrightarrow x$.
So assume that there is $w \in Z$ such that $y \in Z(w)$. As $y
\longrightarrow x \longrightarrow v$, $w \neq v$. Moreover, if $w =
x$, then $\{x, y\}$ is an interval of $T(Z \cup \{y\})$. Since $\{v,
y\}$ is an interval of $T((Z \cup \{y\})-\{x\})$, we would obtain
that $\{x, y, v\}$ is an interval of $T(Z \cup \{y\})$ so that $\{x,
v\}$ would be an interval of $T(Z)$. Therefore, $w \notin \{v, x\}$
and hence $\{v, w\}$ is an interval of $T(Z) - x$. As  $x \in Y(0)$,
it follows from Remark~\ref{r10} that $\{v, w\} = \{n, n+1\}$ so that $v=
n$ and $n \longrightarrow y$. By considering the automorphism
$\sigma^{n+1}$ of $T(Y)$ defined by $ \sigma^{n+1}(i)= i+n+1$, we
obtain that $y\in Y(0)$ and $x\in Y(n+1)$. By considering
$T^{\star}$ instead of $T$, we get $y \in Y(0)$ and $x\in Y(n)$
because the permutation $\pi: i \mapsto -i$ is an isomorphism from
$T(Y)$ onto $T(Y)^{\star}$ by Remark~\ref{r10}. Lastly, by interchanging
$x$ and $y$ in the foregoing, we obtain $n \longrightarrow x$ in
$T^{\star}$ which  means that initially $x\longrightarrow 0$ in $T$.
It follows that the function $Y\cup\{x,
y\}\longrightarrow\{0,\ldots, 2n+2\}$, defined by $ x \mapsto 2n+2$,
$ y \mapsto n+1$, $ i \mapsto i$ for $0\leq i\leq n$ and $ i \mapsto
i+1$ for $n+1\leq i\leq 2n$, realizes  an isomorphism from $T(Y \cup
\{x, y\})$ onto $T_{2n+3}$. Consequently, $T(Y \cup \{x, y\}) \simeq
T_{2n+3}$, with $Y \cup \{x, y\} \subseteq X$, which contradicts the
maximality of $Y$.

Second, suppose that $x\in [Y]$. By interchanging $T$ and
$T^{\star}$, assume that $y \longrightarrow x \longrightarrow Y$.
Consider $Z = (Y - \{v\}) \cup \{y\}$. We have $T(Z) \simeq
T_{2n+1}$ and, by the previous observation, either $x\in [Z]$ or
there is $w\in Z$ such that $x\in Z(w)$. The first instance is not
possible because $y \longrightarrow x \longrightarrow Z-\{y\}$.
Since $y \longrightarrow x \longrightarrow Z-\{y\}$ and hence
$s_{T(Z \cup \{x\})}(x) = 2n$, the second is not possible either.
Indeed, given $w\in Z$, if $x\in Z(w)$, then $s_{T(Z \cup \{x\})}(x)
\in \{n, n+1\}$ because $s_{T(Z)}(w) = n$.

 It follows that $Ext(Y) \cap X \neq
\emptyset$. Set $T(Y)= T_{2n +1}$. By the preceding observation,
 $X = Y \cup \{x\}$, where $x\in Ext(Y) \cap X$. As $\mid X\mid$ is
 even, it follows from Corollary~\ref{c8} that there is $j\in X-\{x\}$ such
 that $T(X)-j$ is indecomposable. By considering the automorphism $\sigma^{2n
 +1-j}$ of $T(Y)$, we can assume that $j=0$. For a contradiction,
 suppose that $T(X) - ~ 0 ~ \simeq ~ T_{2n+1}$. We would have $s_{(T(X) - 0)}(x) =
 n$. Since $s_{(T(Y) - 0)}(i) = n$ for $ 1\leq i \leq n $ and $s_{(T(Y) -
0)}(i) = n-1$ for $n+1\leq i \leq 2n$, we would obtain that
$N^{+}_{(T(X)-0)}(x) = \{1, \ldots, n\} $ so that $\{0, x\}$ would
be a non trivial interval of $T(X)$. Consequently, $T(X)-0$ is not
isomorphic to $T_{2n +1}$. By Proposition~\ref{p4}, a diamond embeds into
$T(X)-0$. It follows from the minimality of $T(X)$ that $T(X)-0$ and
hence $T(Y)-0$ omit $T_{5}$. As $T_{5}$ embeds into $T_{2m +1}-0$
for $m\geq 3$, we get $n=2$.

It remains to verify that  $U_{5}$ and $W_{5}$ embed into $T(X)$.
Since  $x \notin [Y]$, $s_{T(X)}(x) \in \{1,2,3,4\}$. By
interchanging $T$ and $T^{\star}$, assume  that $s_{T(X)}(x) =1$ or
2. First, assume that there is $i\in \mathbb{Z}/5\mathbb{Z}$ such
that  $N^{+}_{T(X)}(x) = \{i\}$. By considering the automorphism
$j\mapsto j-i$ of $T_{5}$, assume that $i=0$. The function
$\mathbb{Z}/5\mathbb{Z}\longrightarrow X-\{3\}$, which fixes $0$,
$1$, $2$, $4$ and which maps $3$ to $x$, is an isomorphism from
$U_{5}$ onto $T(X)-3$. Furthermore, the function
$\mathbb{Z}/5\mathbb{Z}\longrightarrow X-\{2\}$, defined by
$0\mapsto 3$, $1\mapsto 4$, $2\mapsto x$, $3\mapsto 0$ and $4\mapsto
1$, is an isomorphism from $W_{5}$ onto $T(X)-2$. Finally, assume
that there is $i\in \mathbb{Z}/5\mathbb{Z}$ such that
$N^{+}_{T(X)}(x) = \{i, i+1\}$ or $\{i, i+2\}$. If $N^{+}_{T(X)}(x)
= \{i, i+1\}$, then $\{i-1, x\}$ would be an interval of $T(X)$. So,
by considering the automorphism $k\mapsto k-i$ of $T_{5}$, assume
that $N^{+}_{T(X)}(x) = \{0, 2\}$. The function
$\mathbb{Z}/5\mathbb{Z}\longrightarrow X-\{0\}$, defined by
$0\mapsto 2$, $1\mapsto 3$, $2\mapsto 4$,  $3\mapsto x$ and
$4\mapsto 1$, is an isomorphism from $U_{5}$ onto $T(X)-0$.
Furthermore, the function $\mathbb{Z}/5\mathbb{Z}\longrightarrow
X-\{2\}$, defined by $0\mapsto 3$, $1\mapsto 4$, $2\mapsto x$,
$3\mapsto 0$ and $4\mapsto 1$, is an isomorphism from $W_{5}$ onto
$T(X)-2$.
\epr

\section{A new characterization of the critical tournaments}

In this section we   discuss some other questions concerning the
indecomposable subtournaments on 5 and 7 vertices of an
indecomposable tournament. In particular, we obtain a new
characterization of the critical tournaments. In that order, we
recall the following two results concerning the critical
tournaments.

\begin{lem}[\cite{ST}]\label{l11} The indecomposable subtournaments  of $T_{2n+1}$ on at least $5$ vertices, where $n \geq 2$,
are isomorphic to $T_{2m+1}$, where $2 \leq m\leq n$. The same holds
for the indecomposable subtournamants of $U_{2n+1}$ and of
$W_{2n+1}$.

\end{lem}

\begin{lem}[\cite{Y}]\label{l12} Given an indecomposable tournament $T$
of cardinality $\geq5$, $T$ is critical if and only if $T$ omits
any indecomposable tournament on six vertices.\end{lem}

 Let $T$ be an indecomposable tournament of cardinality $\geq
5$. We denote by $I_{5}(T)$ the set of the  elements  of $\{T_{5},
U_{5}, W_{5}\}$ embedding in $T$. By Remark~\ref{r9}, $I_{5}(T) \neq
\emptyset$. By Theorem~\ref{t3}, $I_{5}(T) \neq \{T_{5}, U_{5}\}$ and
$I_{5}(T) \neq \{T_{5}, W_{5}\}$. We characterize the
indecomposable tournaments $T$ such that $I_{5}(T) = \{T_{5}\}$
(resp. $I_{5}(T) = \{U_{5}\}$). The following remark completes
this discussion.

\begin{rem}\label{r13} For $J = \{W_{5}\}$, $\{U_{5}, W_{5}\}$ or $\{T_{5}, U_{5},
W_{5}\}$ and for $n \geq 6$, there exists an indecomposable
tournament $T$ of cardinality $n$ such that $I_{5}(T) = J$.
\\
For $n \geq 5$, the tournaments $E_{n+1}$, $F_{n+1}$ and $G_{n+1}$
defined below on  $\{0,\ldots,$ $n\}$ are indecomposable and satisfy
$I_{5}(E_{n+1}) = \{T_{5}, U_{5}, W_{5}\}$, $I_{5}(F_{n+1})
=\{W_{5}\}$ and $I_{5}(G_{n+1}) =\{U_{5}, W_{5}\}$.
\begin{itemize}
\item $E_{n+1}(\{0, \ldots, 4\}) = T_{5}$ and, for all $5 \leq k \leq n$,
$N_{E_{n+1}(\{0, \ldots, k\})}^{+}(k) = \{k-1\}$;
\item $A(F_{n+1}) = \{(i,j):i+1 < j\ or\ i = j+1\}$;
\item $G_{n}(\{0, \ldots, n-1 \}) = F_{n}$ and $N^{+}_{G_{n+1}}(n) = \{0\}$.

\end{itemize}

\end{rem}

The following is an easy consequence of Theorem~\ref{t2} and of Lemma~\ref{l11}.

\begin{cor}\label{c14} The next two assertions are satisfied by any
indecomposable tournament $T$ of cardinality $\geq 5$.

\begin{enumerate}

\item $T$ is isomorphic to $T_{2n+1}$ for some $n\geq 2$ if and only
if the indecomposable subtournaments of $T$ on 5 vertices are
isomorphic to $T_{5}$.

\item $T$ is isomorphic to $B_{6}$, $P_{7}$ or to $U_{2n+1}$ for
some $n\geq 2$ if and only if the indecomposable subtournaments of
$T$ on 5 vertices are isomorphic to $U_{5}$.

\end{enumerate}

\end{cor}

For all $n \geq 6$, the tournament $F_{n}$ defined in Remark~\ref{r13} is
an indecomposable non critical tournament all the indecomposable
subtournaments of which are isomorphic to $W_{5}$. This leads us to
the following characterization of the tournaments $W_{2n+1}$ and to
the problem  below.

\begin{prop}\label{p15}
Given an indecomposable tournament $T$ of cardinality $\geq 7$, $T$
is isomorphic to  $W_{2n+1}$ for some $n \geq 3$ if and only if the
indecomposable subtournaments on 7 vertices of $T$ are isomorphic to
$W_{7}$.
\end{prop}

\bpr
By Lemma~\ref{l11}, if $T \simeq W_{2n+1}$, where $n \geq 3$, then the
indecomposable subtournaments of $T$ on 7 vertices  are isomorphic
to $W_{7}$. Conversely, assume that  the indecomposable
subtournaments of $T$ on 7 vertices  are isomorphic to $W_{7}$. By
Lemma~\ref{l11}, it suffices to show that $T$ is critical. Clearly,  if
$\mid \!V(T)\! \mid = 7$, then $T \simeq W_{7}$. So assume that
$\mid \!V(T) \!\mid \geq 8$. For a contradiction, suppose that $T$
is not critical. It follows from Lemma~\ref{l12} that there exists
$X\subset V(T)$ such that $\mid \!X\! \mid = 6$ and $T(X)$ is
indecomposable. By Proposition~\ref{p7}, there is $Y\subseteq V(T)$ such
that $X\subset Y$, $\mid \!Y\! \mid = 8$ and $T(Y)$ is
indecomposable. As $\mid \!Y\! \mid$ is even, $T(Y)$ is not
critical. Consider $x \in Y$ such that $T(Y) - x$ is indecomposable.
We have  $ T(Y) - x \simeq W_{7}$ and hence we can denote the
elements of $Y$ by $0,\ldots, 7$ in such a way that  $x = 7$ and
$T(Y)-7 = W_{7}$. By Corollary~\ref{c8}, there is $y \in \{0,\ldots, 6\}$
such that $T(Y)- y$ is indecomposable and thus $T(Y)- y\simeq
W_{7}$. To obtain a contradiction, we verify that $\{y, 7\}$ would
be a non trivial interval of $T(Y)$. By interchanging $T$ and
$T^{\star}$, we can assume that $y\in \{0, 1, 2\}\cup \{6\}$ because
the permutation of $\mathbb{Z}/7\mathbb{Z}$, which fixes $6$ and
which exchanges $i$ and $5- i$ for $0\leq i\leq 5$, is an
isomorphism from $W_{7}$ onto its dual. First, assume that $y= 6$.
We have $ T(Y)-\{6, 7\}= 0<\cdots< 5$. Since $\{1,\ldots, 5\}\cup
\{7\}$ is not an interval of $T(Y)-6$, $7\longrightarrow 0$. As
$\{i, i+1\}$ is not an interval of $T(Y)-6$ for $0\leq i\leq 4$, we
obtain successively that $1\longrightarrow 7$, $7\longrightarrow 2$,
$3\longrightarrow 7$, $7\longrightarrow 4$ and $5\longrightarrow 7$.
Second, assume that $y\in \{0, 1, 2\}$. For $z\in \{0,\ldots,
7\}-\{y, 6\}$, $C_{3}$ embeds into $T(Y)-\{y, z\}$ because $T(\{2i,
2i+1, 6\})\simeq C_{3}$ for $i\in  \{0, 1, 2\}$. It follows that the
isomorphism from $ W_{7}$ onto $T(Y)-y$ fixes $6$. Consequently,
$T(Y)-\{y, 6\}$ is transitive. We have only to check that $T(Y)-\{y,
6\}$ is obtained from the usual total order on $\{0, \ldots, 5\}$ by
replacing $y$ by $7$. If $y =0$, then $7\longrightarrow 1$ because
$1\longrightarrow \{2,\ldots, 6\}$. Thus $T(Y)-\{y, 6\}= 7< 1<\cdots
< 5$. If $y =1$ or $2$, then $\{y-1, y+1\}$ is  an interval of
$T(Y)-\{y, 7\}$. Therefore, $\{y-1, y+1\}$ is not an interval of
$T(\{y-1, y+1, 7\})$ and hence $T(Y)-\{y, 6\}= \cdots < y-1 < 7 <
y+1<\cdots< 5 $.

\epr

From Corollary~\ref{c14} and Proposition~\ref{p15}, we obtain the following
recognition   of the critical tournaments from their indecomposable
subtournaments on 7 vertices.

\begin{cor}\label{c16}

Given an indecomposable tournament $T$, with $\mid\!V(T)\!\mid \geq
7$, $T$ is critical if and only if the indecomposable subtournaments on 7 vertices of $T$ are isomorphic to one and only one of the tournaments $T_{7}$, $U_{7}$ and  $W_{7}$.
\end{cor}

\begin{prob}\label{p17}
Characterize  the indecomposable tournaments all of whose
indecomposable subtournaments on 5 vertices  are isomorphic to
$W_{5}$.
\end{prob}

\end{document}